\documentclass[11pt]{amsart}
\usepackage{amsxtra}
\usepackage{amssymb}
\usepackage{graphicx}
\usepackage{mathrsfs} 
\usepackage{enumerate}
\usepackage{array,color}
\usepackage{hyperref}
\usepackage[all]{xy}
\usepackage{comment}
\theoremstyle{plain}
\newcommand{\cleqn}{\setcounter{equation}{0}}
\newcommand{\clth}{\setcounter{theorem}{0}}
\newcommand {\sectionnew}[1]{\section{#1}\cleqn\clth}
\newtheorem{theorem}{Theorem}[section]
\newtheorem{lemma}[theorem]{Lemma}
\newtheorem{definition-theorem}[theorem]{Definition-Theorem}
\newtheorem{proposition}[theorem]{Proposition}
\newtheorem{corollary}[theorem]{Corollary}
\newtheorem{definition}[theorem]{Definition}
\newtheorem{example}[theorem]{Example}

\newtheorem{remark}[theorem]{Remark}
\newtheorem{conjecture}[theorem]{Conjecture}
\newtheorem{notation}[theorem]{Notation}

\newcommand \bth[1] { \begin{theorem}\label{t#1} }
\newcommand \ble[1] { \begin{lemma}\label{l#1} }

\newcommand \bpr[1] { \begin{proposition}\label{p#1} }
\newcommand \bco[1] { \begin{corollary}\label{c#1} }
\newcommand \bde[1] { \begin{definition}\label{d#1}\rm }
\newcommand \bex[1] { \begin{example}\label{e#1}\rm }
\newcommand \bre[1] { \begin{remark}\label{r#1}\rm }
\newcommand \bcj[1] { \begin{conjecture}\label{j#1}\rm }

\newcommand \bnota[1] { \begin{notation}\label{n#1}\rm }
\renewcommand {\eth} { \end{theorem} }
\newcommand {\ele} { \end{lemma} }

\newcommand {\epr} { \end{proposition} }
\newcommand {\eco} { \end{corollary} }
\newcommand {\ede} { \end{definition} }
\newcommand {\eex} { \end{example} }
\newcommand {\ere} { \end{remark} }
\newcommand {\ecj} { \end{conjecture} }

\newcommand {\enota} { \end{notation} }
\newcommand \thref[1]{Theorem \ref{t#1}}

\newcommand \prref[1]{Proposition \ref{p#1}}

\newcommand \deref[1]{Definition \ref{d#1}}
\newcommand \exref[1]{Example \ref{e#1}}

\newcommand{\vocab}[1]{\emph{#1}}

\def \Cset {{\mathbb C}}

\def \wt {\widetilde}
\def \wh {\widehat}

\def \id { {\mathrm{id}} }

\newcommand{\supp}{\text{supp}}
\begin{document}
%%%%%%%%%%%%%%%%%%%%%%%%%%%%%%%%%%%%%%%%%%%%%%%%%%%%%%%%%%%%%%%%%%%%%%%%%%%
%%%%%%%%%%%%%%%%%%%%%%    Title    %%%%%%%%%%%%%%%%%%%%%%%%%%%%%%%%%%%%%%%%
\title[The restricted DFT]
{The restricted Discrete Fourier Transform}
\author[W. Riley Casper]{W. Riley Casper}
\address{
Department of Mathematics \\
California State University Fullerton \\
Fullerton, CA 92831\\
U.S.A.
}
\email{wcasper@fullerton.edu}

\author[Milen Yakimov]{Milen Yakimov}
\address{
Department of Mathematics, 
Northeastern University, 
Boston, MA 02115, 
U.S.A. and International Center for Mathematical Sciences, Institute of Mathematics and Informatics \\
Bulgarian Academy of Sciences \\ 
Acad. G. Bonchev Str., Bl. 8 \\
Sofia 1113, Bulgaria
}
\email{m.yakimov@northeastern.edu}
\date{}
%\keywords{}
%\subjclass[2020]{Primary ; Secondary }
\begin{abstract} We investigate the restriction of the discrete Fourier transform $F_N : L^2(\mathbb{Z}/N \mathbb{Z}) \to L^2(\mathbb{Z}/N \mathbb{Z})$ to the space $\mathcal C_a$ of functions with support on the discrete interval $[-a,a]$, whose transforms are supported inside the same interval. A periodically tridiagonal matrix $J$ on $L^2(\mathbb{Z}/N \mathbb{Z})$ is constructed having the three properties that it commutes with $F_N$, has eigenspaces of dimensions 1 and 2 only, and the span of its eigenspaces of dimension 1 is precisely $\mathcal C_a$.
The simple eigenspaces of $J$ provide an orthonormal eigenbasis of the restriction of $F_N$ to $\mathcal C_a$. The dimension 2 eigenspaces of $J$ have canonical basis elements supported on $[-a,a]$ and its complement.  These bases give an interpolation formula for reconstructing $f(x)\in L^2(\mathbb{Z}/N\mathbb{Z})$ from the values of $f(x)$ and $\widehat f(x)$ on $[-a,a]$, i.e., an explicit Fourier uniqueness pair interpolation formula. The coefficients of the interpolation formula are expressed in terms of theta functions.  
The collections of simple eigenvalues of $J$ are proved to be strictly greater than the double eigenvalues.
Lastly, we construct an explicit basis of $\mathcal C_a$ having extremal support and leverage it to obtain explicit formulas for eigenfunctions of $F_N$ in $\mathcal C_a$ when $\dim \mathcal C_a \leq 4$.
\end{abstract}
\maketitle

\section{Introduction}
\subsection{The restricted discrete Fourier transform}
The (non-normalized) \vocab{discrete Fourier transform} on 
$\mathbb{Z}/N\mathbb{Z}$ is the linear transformation acting on the Hilbert space $L^2(\mathbb{Z}/N\mathbb{Z})$ of complex-valued functions $f(k)$ on $\mathbb{Z}/N\mathbb{Z}$ given by
$$F_N:f(k)\mapsto \wh f(k) := \sum_{j=0}^{N-1}e^{-2\pi ijk/N}f(j).$$
The eigenfunctions of the discrete Fourier transform play an important role in defining fractional Fourier transforms and are connected with theta functions.
Each eigenvalue of $F_N$ has a large multiplicity, leading to many choices for an eigenbasis.
In general what constitutes a nice choice of basis for the eigenvectors of $F_N$ depends on the intended application and several methods have been suggested for obtaining one.
This problem has been studied by many authors \cite{landau,good,mcclellan2,auslander,dickinson,GHS}.
%It was pointed out by Good and McClellan that the 
%The problem of obtaining the eigenvalues of the discrete Fourier transform and their multiplicities
%to be equivalent to a problem originally considered by Gauss.

This paper is dedicated to the study of the restriction of the discrete Fourier transform to the space of functions supported on the discrete interval
\[
[-a, a] := \{-a,1-a,\dots,a-1,a\},
\]
whose transforms are supported inside the same interval. 
Put precisely, this is the restriction of $F_N$ to the space
\[
\mathcal C_a := L^2([-a,a]) \cap F^{-1}_N( L^2([-a,a])), 
\]
where $L^2(X)$ denotes the subspace of functions in $L^2(\mathbb{Z}/N\mathbb{Z})$ with support contained in a subset $X\subseteq \mathbb{Z}/N\mathbb{Z}$.
In stark contrast to the setting of continuous intervals and Fourier transforms on the real line or the circle, the space $\mathcal C_a$ can be nontrivial. We will call the restriction
$$F_N|_{\mathcal C_a}$$
the {\vocab{restricted discrete Fourier transform.}}

\subsection{Spectral analysis of the restricted Fourier transform}
Our analysis of the restricted Fourier transform is based on leveraging the two operators on $L^2(\mathbb{Z}/N\mathbb{Z})$ represented by the 
two periodically tridiagonal matrices
$$
J_0 := \left[\begin{array}{ccccc}
b_0 & 1 &  0  & \dots & 1\\
1 & b_1 & 1 & \dots & 0\\
 0  & 1 & b_2 & \dots & 0\\
\vdots & \vdots & \vdots & \ddots & \vdots\\
1  &  0  &  0  & \dots & b_{N-1}
\end{array}\right]
\quad\text{and}\quad
J_1 := \left[\begin{array}{ccccc}
0 & a_1 &  0  & \dots & a_N\\
a_1 & 0 & a_2 & \dots & 0\\
 0  & a_2 & 0 & \dots & 0\\
\vdots & \vdots & \vdots & \ddots & \vdots\\
a_N  &  0  &  0  & \dots & 0
\end{array}\right],
$$
in the standard basis of $L^2(\mathbb{Z}/N\mathbb{Z})$, where 
$$b_k := 2\cos \frac{2\pi k}{N}\quad\text{and}\quad a_k := \cos \frac{\pi (2k-1)}{N}.$$
Both $J_0$ and $J_1$ commute with $F_N$, and one can show that together they generate the algebra of all matrices commuting with $F_N$. Here, we are concerned with the operator 
\begin{equation}
\label{operator-J}
J := J_1 - \cos\frac{\pi(2a+1)}{N}J_0
\end{equation}
for an integer $a$. 
Studying the spectra of $J$ restricted to the space of functions supported on the discrete interval $[-a,a]$ and its complement
$$
[-a,a]' = \left( \mathbb{Z}/N \mathbb{Z} \right) \backslash [-a,a]
$$
leads to the following set of results.
\smallskip

\noindent
{\bf{Spectral Theorem.}} {\em{
Let $0\leq a \leq (N-1)/2.$ 
The space $\mathcal C_a$ is nontrivial if and only if $a\geq (N-1)/4$. The following hold for $a\geq (N-2)/4$:}}
\begin{enumerate}
\item[(a)] {\em{The dimension of $\mathcal C_a$ is $r=4a+2-N$ and has basis
\begin{equation}
\label{psi-basis}
%\{\psi_a(x),e^{-2\pi ix/N}\psi_a(x),\dots,e^{-2\pi i(r-1)x/N}\psi_a(x)\},
\{e^{-2\pi i jx/N}\psi_a(x): 0\leq j < r\}
\end{equation}
where the function $\psi_a(x)$ has support in $[-a,a]$ with
\[
\psi_a(x) = 
e^{i\pi (r-1)x/N}\prod_{k=1}^{N-2a-1}\sin\left(\frac{\pi (a+k-x)}{N}\right).
\]}}
\item[(b)] {\em{The operator $J$ preserves the subspaces $L^2([-a,a])$ and $L^2([-a,a]')$ of $L^2(\mathbb{Z}/N \mathbb{Z})$ and its restrictions to those subspaces 
are given by Jacobi matrices in the standard bases. In particular, the spectra of 
$$J|_{L^2([-a,a])} \quad \mbox{and} \quad J|_{L^2([-a,a]')}$$  
are simple.}} 
\item[(c)] {\em{Each eigenvalue of $J|_{L^2([-a,a]')}$ is an eigenvalue of $J|_{L^2([-a,a])}$.}} 
\end{enumerate}
{\em{By parts (b)--(c), there is a unique (up to rescaling by $\pm 1$) choice for a real, orthonormal eigenbasis for $J$ of the form
\begin{equation}
\label{rhovarphi}
\{\rho_1(x),\dots,\rho_r(x),\varphi_1(x),\wt\varphi_1(x),\dots,\varphi_s(x),\wt\varphi_s(x)\},
\end{equation}
where $s := (N - r)/2 = N-2a-1$, 
$$J\rho_k(x) = \mu_k\rho_k(x), \; \; J\varphi_j(x) = \nu_j\varphi_j(x), \;\; \text{and} \;\; J\wt\varphi_j(x) = \nu_j \wt\varphi_j(x), \;\;
1 \leq k \leq r, 1, \leq j \leq s,$$
the eigevalues $\{\mu_1, \ldots, \mu_r, \nu_1, \ldots, \nu_s\}$ are distinct, the functions $\rho_k(x)$ and $\varphi_j(x)$ are supported in $[-a,a]$ and $\wt\varphi_j(x)$ are supported in the complement $[-a,a]'$. The sets 
$$
\{ \rho_1(x), \ldots, \rho_r(x), \varphi_1(x), \ldots, \varphi_s(x) \} 
\quad \mbox{and} \quad
\{\wt\varphi_1(x),\dots,\wt\varphi_s(x) \}
$$
are complete collections of orthonormal eigenvectors of the restrictions of $J$ to $L^2([-a,a])$ and $L^2([-a,a]')$, respectively.}}
\begin{enumerate} 
\item[(d)] {\em{The operator $J$ preserves $\mathcal C_a \subseteq L^2([-a,a])$ and has simple spectrum on it. The set
$$
\{\rho_1(x), \ldots, \rho_r(x) \}
$$
is a joint orthonormal eigenbasis of $\mathcal C_a$ for the commuting operators 
$F_N|_{\mathcal C_a}$ and $J|_{\mathcal C_a}$.}} 
\item[(e)] The collections of double eigenvalues 
$\{\nu_1, \ldots, \nu_s\}$ and simple eigenvalues $\{\mu_1, \ldots, \mu_r\}$ are strictly separated: 
\[
\max_{1\leq j\leq s}\nu_j \leq 2\cos\left(\frac{\pi}{N}\right)+4\cos\left(\frac{\pi (2a+1)}{N}\right) \leq \min_{1\leq k\leq r}\mu_k.
\]
\item[(f)] {\em{The eigenvalues and multiplicities of the restriction of $F_N|_{\mathcal C_a}$ are given by the entries in the table in Figure \ref{eigenvalues of F}.}}
\end{enumerate}
\bigskip

\begin{figure}[htp]
    \centering
    \begin{tabular}{|c|c|c|c|c|}\hline
       $N$     & $\lambda=\sqrt{N}$ & $\lambda=-i\sqrt{N}$ & $\lambda=-\sqrt{N}$ & $\lambda=i\sqrt{N}$\\
       \hline
       $4m-2$  &     $a-m+1$  &     $a-m+1$   &     $a-m+1$   &     $a-m+1$    \\
       $4m-1$  &     $a-m+1$  &     $a-m+1$   &     $a-m+1$   &     $a-m$      \\
       $4m$    &     $a-m+1$  &     $a-m+1$   &     $a-m$     &     $a-m$      \\
       $4m+1$  &     $a-m+1$  &     $a-m$     &     $a-m$     &     $a-m$      \\
       \hline
    \end{tabular}
    \caption{The eigenvalue multiplicities of the restriction of the discrete Fourier transform to $\mathcal C_a$ for $a\geq m$.}
    \label{eigenvalues of F}
\end{figure}
We note that the case of $a = (N-2)/4$ results in a trivial space $\mathcal C_a$ but nontrivial statement in part (c) and the results that follow. 

The problem of obtaining the multiplicities of the eigenvalues of the discrete Fourier transform has been studied by many authors \cite{landau,good,mcclellan2,auslander,dickinson} and was pointed out by Good and McClellan to be equivalent to a problem originally considered by Gauss \cite{mcclellan}.
The multiplicities of the eigenvalues of the restricted discrete Fourier transform obtained in part (f) of the Spectral Theorem generalize these previously known results.

The eigenfunctions of the operator $J$ are examples of prolate spheroidal wave functions of the discrete Fourier transform. The discrete continuous case (Fourier series) was studied by Slepian \cite{slepian,slepian2}.
The eigenfunctions of $J$ belonging to $L^2([-a,a])$ are examples of prolate spheroidal wave functions for the \emph{finite} Fourier transform which limit asymptotically to the prolates studied by Slepian.  These have been explored under various names (such as discrete-discrete prolates and periodic discrete prolates) by several authors \cite{grunbaum3,jain,wilson,pei} for connections to applied settings, which are very different from our methods. Among many differences, we have explicit expressions for eigenfunctions in certain cases and applications to interpolation formulas linked to theta functions and Fourier uniqueness pairs, which we present next. All previous works studied the restrictions of $J$ to $L^2([-a,a])$ and $L^2([-a,a]')$ in isolation, while we investigate the interrelation between the two restrictions and show that it governs the restricted discrete Fourier transform.

The idea of using an operator which commutes with $F_N$ to find some eigenfunctions of the discrete Fourier transform was used by Gr\"{u}nbaum \cite{grunbaum} and Dickinson and Steiglitz \cite{dickinson}. However, it was not realized that a complete spectral analysis of the restricted discrete Fourier Transform can be obtained from the simple eigenspaces of a commuting difference operator. 

\bre{different bases} The two bases 
$$\{\psi_a(x),e^{-2\pi ix/N}\psi_a(x),\dots,e^{-2\pi i(r-1)x/N}\psi_a(x)\} \quad \mbox{and} \quad \{\rho_1(x), \ldots, \rho_r(x) \}$$
of $\mathcal C_a$ from parts (a) and (d) of the Spectral Theorem are of very different nature. The second is orthonormal and consists 
of the eigenvectors of a Jacobi matrix. The first is not orthonormal and consists of functions $f(x)$ of extremal support on $\mathbb{Z}/N \mathbb{Z}$ 
in the sense that for all functions $g(x)$ on $\mathbb{Z}/N\mathbb{Z}$, 
$$
\supp(g)\subseteq \supp(f) \quad \mbox{and} \quad \supp(\hat g)\subseteq\supp(\hat f)
\quad \Rightarrow \quad g(x) = \mu f(x) \; \; \mbox{for some} \; \; \mu \in \Cset.
$$
The last property is proved in \thref{extremal-support-psi}. 
\ere
\medskip

\noindent
{\bf{Extremal Cases for $\dim \mathcal C_a$:}}
\begin{enumerate}
\item If $N$ is odd, choosing $a = (N-1)/2$ gives 
$$
[-a,a] = \mathbb{Z} / N \mathbb{Z},
\quad 
\mbox{and thus}, 
\quad
\mathcal C_a = L^2(\mathbb{Z}/N\mathbb{Z}).
$$
In this case $r=N$, $s=0$ in the Spectral Theorem and 
$$
\{\rho_1(x), \ldots, \rho_N(x) \}
$$
is a complete orthonormal collection of eigenvectors of both $F_N$ and $J$.
\item The basis \eqref{psi-basis} of $\mathcal C_a$ from part (a) of the Spectral Theorem is not an eigenbasis of $F_N|_{\mathcal C_a}$.
However, one can easily obtain explicit eigenbases for small values of $\dim \mathcal C_a$. This is done in Section \ref{sec:low dim} when
$$
1 \leq \dim \mathcal C_a \leq 4. 
$$ 
It seems that among these eigenfunctions only 2 were found before: Kong \cite{kong} found a generator of $\mathcal C_a$ in the case when $\dim \mathcal C_a =1$ and one of the two eigenfunctions when $\dim \mathcal C_a =2$. 
\end{enumerate}

%this was never elevated to obtaining 
%in both papers multiplicity of some of the eigenvectors in the commuting operator prevent the strategy from obtaining a complete eigenbasis of a restriction of the 
\subsection{Interpolation and Fourier uniqueness pairs}
One particular consequence of the previous theorem is that a pair of identical discrete intervals of the form $[-a,a]$ for $a\geq (N-2)/4$ forms a \vocab{Fourier uniqueness pair} for the group 
$\mathbb{Z}/N \mathbb{Z}$.
By this we mean a pair of sets $A,B\subseteq\mathbb{Z}/N \mathbb{Z}$ where knowledge of $f(x)$ on $A$ and $\wh f(x)$ on $B$ determines the entire function $f(x)$.
Fourier uniqueness pairs for discrete subsets of the real line were introduced by Cohn, Kumar, Miller, Radchenko, and Viazovska in the context of sphere packing problems in dimension $8$ and $24$ \cite{Viazovska,Cohn}.
Recently, they have also been connected with $L$-functions and the Riemann hypothesis \cite{Bondarenko}.
In this context, the \vocab{Fourier interpolation problem} is the problem of determining $f(x)$ from the known data on $A$ and $B$ in terms of an expansion of certain ``magic" functions. 

The eigenfunctions of $J$ have a natural immediate application to the Fourier interpolation problem as the desired magic functions as shown in our second main theorem.
For its statement, we introduce some additional notation. 
Denote the projection
$$
P_a : L^2(\mathbb{Z}/N\mathbb{Z}) \to L^2([-a,a]).
$$
By a slight abuse of terminology, we will call the operator 
\begin{equation}
\label{tblimF}
F_N^a := P_a F_N P_a : L^2(\mathbb{Z}/N\mathbb{Z}) \to L^2([-a,a])
\end{equation}
the {\em{time-band limited discrete Fourier transform}}.
The operator 
\begin{equation}
\label{tblimF2}
B_N^a := (F_N^a)^*(F_N^a) = P_a F_N^* P_a F_N P_a: L^2(\mathbb{Z}/N\mathbb{Z}) \to L^2([-a,a])
\end{equation}
is called the {\em{time-band limiting operator}}.

By way of definition, 
\begin{equation}
\label{tbF-F}
B_N^a|_{\mathcal C_a} = N\cdot\text{id}_{\mathcal C_a}
\quad\text{and}\quad
F_N^a|_{\mathcal C_a} = F_N|_{\mathcal C_a}.
\end{equation}
\medskip

\noindent
{\bf{Interpolation Theorem.}} {\em{
Let $(N-2)/4 \leq a \leq (N-1)/2$ be an integer and set $r = 4a+2-N$ and $s = N-2a-1$. Then in the notation of the Spectral Theorem, we have the following:}} 
\begin{enumerate}
\item[(a)] {\em{For all $1 \leq j \leq s$, 
$$\left[\begin{array}{c}
F_N \varphi_j(x)\\ F_N\wt\varphi_j(x)
\end{array}\right] = 
\left[\begin{array}{cc}
\alpha_j & \beta_j\\
\beta_j & -\alpha_j
\end{array}\right]
\left[\begin{array}{c}
 \varphi_j(x)\\ \wt\varphi_j(x)
\end{array}\right]
$$
for some nonzero complex numbers $\alpha_j$ and $\beta_j$. They are either both real or both imaginary and $|\alpha_j|^2+|\beta_j|^2=N$.}}
\item[(b)] {\em{The time-band limited Fourier transform $F_N^a$ and the time-band limiting operator $B_N^a$ given by \eqref{tblimF}-\eqref{tblimF2} both commute with $J$. The set \eqref{rhovarphi} is a joint orthonormal eigenbasis for $J$, $F_N^a$ and $B_N^a$ acting on $L^2(\mathbb{Z}/N \mathbb{Z})$:
\begin{itemize}
\item $\rho_k(x)$, $1 \leq k \leq r$ are eigenfunctions of $F_N^a$ with the same eigenvalues as $F_N$, and eigenfunctions of $B_N^a$ with eigenvalue $N$. 
\item $F_N^a(\varphi_j(x)) = \alpha_j \varphi_j(x)$ and $F_N^a \wt{\varphi}_j(x) = 0$ for all $1 \leq j \leq s$.
\item $B_N^a(\varphi_j(x)) = \lvert\alpha_j\rvert^2 \varphi_j(x)$ and $B_N^a \wt{\varphi}_j(x) = 0$ for all $1 \leq j \leq s$.
\end{itemize}
}}
\item[(c)]
{\em{For any $f(x)\in L^2(\mathbb{Z}/N \mathbb{Z})$, we can write
%\begin{align}\label{eqn:interpolation}
%f(x)  &= \sum_{k=1}^s\frac{1}{\beta_k}\left(\sum_{y\in [-a,a]} (\wh f(y) - \alpha_k f(y))\varphi_k(y)\right)\wt \varphi_k(x),\quad\text{for all}\ x\notin[-a,a].
%\end{align}
$$
f(x) = \sum_{y\in [-a,a]}\left(v_y(x)f(y) + w_y(x)\wh f(y)\right),\quad\text{for all}\ x\notin[-a,a]
$$
for the functions $v_{-a}(x),\dots,v_a(x)$ and $w_{-a}(x),\dots,w_a(x)$ defined by
$$
v_y(x)=\sum_{j=1}^s\frac{-\alpha_j}{\beta_j}\varphi_j(y)\wt\varphi_j(x)
\quad\text{and}\quad
w_y(x)=\sum_{j=1}^s\frac{1}{\beta_j}\varphi_j(y)\wt\varphi_j(x).
$$
}}
\item[(d)] {\em{The functions $v_y(x)$ and $w_y(x)$ in the interpolation formula can be expressed in terms of Wronskians of the Jacobi theta function $\vartheta(z,\tau)$.
Specifically, for $\theta(x,\tau) := e^{i\pi \tau x^2/N}\vartheta(x\tau,N\tau)$, we have
$$
v_y(x) = \frac{W\left(\theta(-a,\tau),\dots \theta(x,\tau),\dots,\theta(a,\tau),\vartheta(-a/N,\tau/N),\dots,\vartheta(a/N,\tau/N)\right)}{W\left(\theta(-a,\tau),\dots,\theta(a,\tau),\vartheta(-a/N,\tau/N),\dots,\vartheta(a/N,\tau/N)\right)},
$$
where $\theta(x,\tau)$ is occurring in the $(y+a+1)$'th position, and 
$$
w_y(x) = \frac{W\left(\theta(-a,\tau),\dots,\theta(a,\tau),\vartheta(-a/N,\tau/N),\dots,\theta(x,\tau),\dots,\vartheta(a/N,\tau/N)\right)}{W\left(\theta(-a,\tau),\dots,\theta(a,\tau),\vartheta(-a/N,\tau/N),\dots,\vartheta(a/N,\tau/N)\right)},
$$
where $\theta(x,\tau)$ is occurring in the $(y+3a+2)$'th position.}}
\end{enumerate}
\bigskip

In the Spectral Theorem we saw that the eigenfunctions $\{\rho_k(x)\}_{k=1}^r$ of $J$ recover an eigenbasis of $F_N|_{\mathcal C_a}$. Here we see that 
the rest of the eigenfunctions $\{\varphi_j(x),\wt\varphi_j(x)\}_{j=1}^s$ of $J$ play a central role in the Fourier interpolation problem.

These expressions can be used to derive interesting relationships between theta functions and their derivatives as illustrated in \exref{N2}. 
This again highlights the unique utility of the eigenfunctions of the matrix $J$.

\bre{FN-eigenbasis}
The combination of part (d) of the Spectral Theorem and part (a) of the Interpolation Theorem also gives that
$$\Big\{ \rho_k(x), 
\left(\alpha_j\pm \sqrt{\alpha_j^2+\beta_j^2} \; \right)\varphi_j(x) + \beta_j\wt\varphi_j(x)
: 1 \leq k \leq r, 1 \leq j \leq s 
\Big\}
$$
is an eigenbasis of $F_N$ (acting on 
$L^2(\mathbb{Z}/ N \mathbb{Z})$).
\ere

\noindent 
{\bf{Acknowledgements.}} The research of W.R.C. has been supported by an AMS-Simons Research Enhancement Grant and RSCA intramural grant 0359121 from CSUF, and that of M.Y. by the Bulgarian Science Fund grant KP-06-N92/5 and the Ministry of Education and Science grant DO1-239/10.12.2024, the Simons Foundation grant SFI-MPS-T-Institutes-00007697 and NSF grant DMS–2200762.

\sectionnew{Spectra of the restricted Fourier transform}

In this section we prove all statements in the Spectral Theorem except the second statement in part (a) of it.

\bpr{Ca} Let $0\leq a \leq (N-1)/2.$ The space $\mathcal C_a$ is nontrivial if and only if $a\geq (N-1)/4$, in which case it has dimension
\begin{equation}
\label{dim-Ca}
\dim \mathcal C_a =4a+2-N.
\end{equation}
\epr
\begin{proof} Consider a function $f(x) \in L^2(\mathbb{Z}/N \mathbb{Z})$ supported on $[-a,a]$. Its Fourier transform will be supported on $[-a,a]$
if and only if the vector $[f(-a), \ldots, f(a)]^t$ is in the kernel of the $(N-2a-1)\times (2a+1)$ matrix
$$
\Big[ e^{-2\pi ijk/N} \Big]_{j \notin [-a,a], k \in [-a,a]}. 
$$
This matrix has full rank since each of its principal submatrices is the product of a Vandermonde matrix and a nondegenerate diagonal matrix. 
%or the transpose of such a product.  
If $a < (N-1)/4$, then $N-2a -1 \geq 2a +1$ and the kernel of the matrix is trivial. If $a \geq (N-1)/4$ then the kernel of the matrix 
has dimension 
$$
2 a + 1 - (N - 2a -1) = 4a +2 - N.
$$ 
\end{proof}
The proposition gives the first statement in the Spectral Theorem and the first statement in part (a) of that theorem. 

Denote the shift operator on $L^2(\mathbb{Z}/ N \mathbb{Z})$:
$$\delta_x^n\cdot f(x) = f(x+n).$$
The operator \eqref{operator-J} is given by 
\begin{align*}
J = A(x)\delta_x + B(x) + A(x-1)\delta_x^{-1},
\end{align*}
where the coefficient functions are
\begin{align*}
  A(x) &= \cos\left(\frac{\pi (2x+1)}{N}\right)-\cos\left(\frac{\pi (2a+1)}{N}\right),\\
  B(x) &= - 2\cos\left(\frac{\pi (2a+1)}{N}\right)\cos\left(\frac{2\pi x}{N}\right).
\end{align*}
In particular, 
$$A(a) = 0 \quad \mbox{and} \quad A(-a-1)=0.$$
Therefore, for a function $f(x)\in L^2(\mathbb{Z}/N\mathbb{Z})$ whose support is contained in $[-a,a]$, $(Jf)(x)$ has support contained in the same set.
Thus $J$ preserves the subspace $L^2([-a,a])$ of $L^2(\mathbb{Z}/N\mathbb{Z})$.
Moreover, since $J$ is selfadjoint, it must also preserve the orthogonal complement
$$L^2([a,a])^\perp = L^2 ( [-a,a]').$$

The restriction of $J$ to $L^2([-a,a])$ is a symmetric, tridiagonal matrix with strictly positive entries for the off-diagonal elements, and therefore has simple spectrum.
Likewise the restriction of $J$ to $L^2([-a,a]')$ has simple spectrum.  This shows part (b) of the Spectral Theorem. 
%Thus if $J$ commutes with $F_N$, then its eigenfunctions will automatically be eigenfunctions of $F_N$.  This is the first thing we prove.
\bpr{JFN-comm}
The operators $J_0$ and $J_1$ commute with $F_N$. In particular, $J$ and $F_N$ commute.
\epr
\begin{proof}
One computes directly that
$$F_N^{-1}\delta_x^{\pm 1} F_N f(x) = e^{\mp 2\pi i x/N}f(x)$$
and
$$F_N^{-1}e^{\pm 2\pi ix/N} F_N f(x) = \delta_x^{\pm 1} f(x).$$
Therefore
$$J_0 = \delta_x + F_N^{-1}\delta_xF_N + F_N^{-2}\delta_xF_N^2 + F_N^{-3}\delta_xF_N^3$$
and
$$J_1 = \frac{1}{2}e^{i\pi/N}\left(e^{2\pi ix}\delta_x + F_N^{-1}e^{2\pi ix/N}\delta_xF_N + F_N^{-2}e^{2\pi ix}\delta_xF_N^2 + F_N^{-3}e^{2\pi ix}\delta_xF_N^3\right).$$
Since $F_N^4$ is the identity matrix, it is clear from this that $J_0$ and $J_1$ commute with $F_N$.  Hence $J$ commutes with $F_N$.
\end{proof}
The formulas for $J_0$ and $J_1$ in terms of a sum of permutations of powers of $F_N$ reveal how the operators were found in the first place.
Given any operator $T$ on $L^2(\mathbb{Z}/N\mathbb{Z})$, the equation
$$T + F_N^{-1}TF_N + F_N^{-2}TF_N^2 + F_N^{-3}TF_N^3$$
defines an operator commuting with $T$.
However, $J_0$ and $J_1$ are in a way even more fundamental.
We can prove that the operators $J_0$ and $J_1$ are complete in the sense that they generate the algebra of all operators on $L^2(\mathbb{Z}/N\mathbb{Z})$ commuting with $F_N$. The proof of this fact will appear elsewhere since the fact does not play a role in this paper.
\bpr{simple spectrum}
The operator $J$ preserves $\mathcal C_a$ and its restriction to $\mathcal C_a$ has simple spectrum.  In particular, because of \prref{JFN-comm}, 
an eigenbasis of $J|_{\mathcal C_a}$ is also an eigenbasis of $F_N|_{\mathcal C_a}$.
\epr
\begin{proof}
The fact that $J$ preserves $\mathcal C_a$ follows immediately from the definition of $\mathcal C_a$ and the facts that $J$ and $F_N$ commute
and $J$ preserves $L^2([-a,a])$. 
Furthermore, since $\mathcal C_a$ is a subspace of the space $L^2([-a,a])$ on which $J$ has simple spectrum, the restriction of $J$ to $\mathcal C_a$ also has simple spectrum.  It follows that the eigenfunctions of $J|_{\mathcal C_a}$ are automatically eigenfunctions of $F_N|_{\mathcal C_a}$.
\end{proof}

However, $J$ itself does not have simple spectrum since there will be overlap between the eigenvalues in each of the restrictions.
To see this, consider the projection operators
$$P_a: L^2(\mathbb Z/N)\rightarrow L^2([-a,a]) \quad \text{and} \quad P_a^\perp = \id-P_a: L^2(\mathbb Z/N)\rightarrow L^2([-a,a]').$$
The next proposition gives part (c) of the Spectral Theorem.

\bpr{repeated eigenvalues of J}
Suppose $a\geq (N-2)/4$.
An eigenvalue $\lambda$ of $J$ has multiplicity greater than $1$ if and only if it is an eigenvalue of the restriction of $J$ to $L^2([-a,a]')$.
In this case, it has multiplicity $2$.
Moreover, if $f\in L^2([-a,a]')$ is an eigenfunction with eigenvalue $\lambda$, then so is $P_a\wh f\in L^2([-a,a])$.
\epr
\begin{proof}
If $a\geq (N-2)/4$, then there does not exist $f\in L^2([-a,a]')$ with $\wh f\in L^2([-a,a]')$. Indeed, for such a function $f(x)$, 
the vector $\big[ f(k) : k \notin [-a,a] \big]^t$ will be in the kernel of the $(2a+1) \times (N-2a -1)$ matrix 
$$
\Big[ e^{-2\pi ijk/N} \Big]_{j \in [-a,a], k \notin [-a,a]}. 
$$
This matrix has trivial kernel since it has full rank and $2a+1 \geq N- 2a -1$ 
(as in \prref{Ca}, each of its principal submatrices is the product of a Vandermonde matrix and a nondegenerate diagonal matrix).
Alternatively, when $N$ is prime it is an automatic consequence of \cite{tao}.
Therefore, if $f\in L^2([-a,a]')$ is a eigenfunction of $J$ with eigenvalue $\lambda$, then $P_a \wh f\in L^2([-a,a])$ is nonzero.
The matrix $J$ commutes with $P_a$ and $F_N$, and thus $P_a\wh f$ is an eigenfunction of $J$ with the same eigenvalue $\lambda$.
The rest of the statement of the proposition follows immediately. 
\end{proof}

If instead $f\in L^2([-a,a])$ is an eigenfunction with eigenvalue $\lambda$ which does not appear as one of the eigenvalues of the restriction of $J$ to $L^2([-a,a]')$, then $P_a^\perp \wh f$ must be zero.  This means that $f\in \mathcal C_a$, and consequently $f$ is an eigenfunction of $F_N$ also.
This proves the following proposition. 
\bpr{dimension of restriction}
Let $(N-2)/4\leq a\leq (N-1)/2$.
Then $\lambda$ is an eigenvalue of the restriction of $J$ to $\mathcal C_a$ if and only if $\lambda$ appears as an eigenvalue of $J$, but not an eigenvalue of $J|_{L^2([-a,a]')}$, 
or equivalently $\lambda$ is an eigenvalue of $J|_{L^2([-a,a])}$ but not of $J_{L^2([-a,a]')}$. 
\epr
The fact that the restrictions $J|_{L^2([-a,a])}$ and $J_{L^2([-a,a]')}$ have simple spectra and \prref{dimension of restriction} give a second proof of the dimension formula \eqref{dim-Ca}:
$$
\dim\mathcal C_a = \dim L^2([-a,a]) - L^2([-a,a]') = (2a+1) - (N- 2a+1) = 4a+2-N.
$$

Part (d) of the Spectral Theorem follows from Propositions \ref{psimple spectrum} and \ref{pdimension of restriction}, and  
the facts that the spectra of $J_{L^2([-a,a])}$ and $J|_{L^2([-a,a]')}$ are simple and the operator $J$ is selfadjoint. 

%\begin{proof}
%The preceding discussion proves everything in the statement of the theorem, except for the value of the dimension.
%To see that, first notice that the number of simple eigenvalues of $J$ is exactly the dimension of $\mathcal C$.
%Moreover, $J$ restricted to $L^2([-a,a])$ will have $\dim(L^2([-a,a])) = 2a+1$ distinct eigenvalues.  Likewise the restriction of $J$ to $L^2([-a,a]')$ will have $N-2a-1$ distinct eigenvalues, all of which are duplicates.
%Therefore the number of unique eigenvalues is $\dim(\mathcal C_a) = 2a+1-(N-2a-1) = 4a+2-N$.
%\end{proof}

%\bre{Vandermonde}
%We can get the value of the dimension of $\mathcal C_a$ easily an entirely different way.
%First decompose the Fourier transform as
%$$F_N = P_aF_NP_a + P_a^\perp F_NP_a + P_aF_N P_a^\perp + P_a^\perp F_NP_a^\perp.$$
%Then $\mathcal C_a = \ker (P_a^\perp F_NP_a)\cap L^2([-a,a])$, which in terms of the matrix representation of the DFT is the kernel of a $(N-2a-1)\times (2a+1)$ Vandermonde matrix.
%Since submatrices of Vandermonde matrices have full rank, the dimension estimate follows immediately.  This is essentially the argument of Tao \cite{tao} for arbitrary subsets when $N=p$ is prime.
%\ere

Next, we determine the multiplicities of the eigenvalues of the restriction of $F_N$ to $\mathcal C_a$, thus proving part (f) of the Spectral Theorem.
The strategy is to use the action of a twisted version of the operator $J_0$, i.e., 
$$J_0^{(\lambda)} := \delta_x+\lambda^{-1} F_N^{-1}\delta_xF_N+\lambda^{-2}F_N^{-2}\delta_xF_N^2+\lambda^{-3}F_N^{-3}\delta_x F_N^3,$$
where $\lambda$ is a fourth root of unity.
This operator is nonzero and satisfies the commutation relation
$$F_N^{-1}J_0^{(\lambda)} F_N = \lambda J_0^{(\lambda)}.$$
Therefore, for any $f\in L^2(\mathbb{Z}/N\mathbb{Z})$
$$\wh{J_0^{(\lambda)} f} = \lambda J_0^{(\lambda)}\wh f.$$

\bpr{eigenvalues of restriction}
The eigenvalues of the restriction of the discrete Fourier transform restricted to $\mathcal C_a$ are given by the values in the table in Figure \ref{eigenvalues of F}.
\epr
\begin{proof}
Let $N=4m+2-d$ with $1 \leq d \leq 4$. 
We proceed by induction on $a$.
The base case of $a=m$ is proved in 
Section \ref{sec:low dim} where we explicitly diagonalize $F_N$ on $\mathcal C_a$ in those 4 cases. In addition we show that in each of those cases, $\mathcal C_a$ contains an eigenfunction with eigenvalue $\sqrt{N}$, which does not vanish at $-a$ and $a$.

As an inductive assumption, assume that $\mathcal C_a$ has the right multiplicities and also contains an eigenfunction $f(x)$ with eigenvalue $\sqrt{N}$, which is nonzero for $x=\pm a$.
Then for $\lambda$ a fourth root of unity, the function
$$f_\lambda(x) := (J_0^{(\lambda)} f)(x)$$
satisfies $f_\lambda(a+1)\neq 0$. 
Therefore it belongs to $\mathcal C_{a+1}$, but not $\mathcal C_a$. Moreover,
$$\wh {f_\lambda}(x) = F_N(J_0^{(\lambda)} f)(x) = \lambda^{-1} J_0^{(\lambda)} F_N f(x) = \lambda^{-1}\sqrt{N} f_\lambda(x).$$
Thus the multiplicity of each eigenvalue increases by at least $1$ in passing from $\mathcal C_a$ to $\mathcal C_{a+1}$.
Since $\dim \mathcal C_{a+1} - \dim \mathcal C_a =4$, this describes the entire change to the spectrum.
Finally, $f_1(x)$ is nonvanishing on $x=\pm a$, so by induction our theorem is true.
\end{proof}

\bre{eigenfunction strategy}
If one is only interested in finding any eigenbasis whatsoever for the restrictions of $F_N$ on $\mathcal C_a$, then the previous proof suggests a simple method based on the repeated application of the operator $J_0^{(\lambda)}$.  However, the resultant basis is undesirable from a numerical standpoint, since the eigenfunctions of $F_N$ generated this way with the same eigenvalue will have a high covariance, making them numerically difficult to tell apart.
In contrast, the spectrum of $J$ will be simple, so the eigenfunctions it generates will be orthogonal.
\ere
\sectionnew{An extremal basis for $\mathcal C_a$}
\label{sec:basis for C}
In this section, we prove part (a) of the Spectral Theorem which amounts to constructing the basis \eqref{psi-basis} of $\mathcal C_a$. We furthermore prove that the elements of this basis have extremal support in the sense of \deref{extremal-support} below. 

Set $N=4m+2-d$ with $1\leq d\leq 4$. As we saw in the previous section, the minimal value of $a$ for which $\mathcal C_a$ is nontrivial is $a=m$, and in that case $\dim \mathcal C_m = d$. Denote a primitive 
$N$'th root of unity
$$\xi:=e^{2\pi i /N}.$$ 
We will use the Gaussian binomial coefficients
$$\binom{n-1}{x}_\xi = \frac{(1-\xi^{n-1})(1-\xi^{n-2})\dots(1-\xi^{n-x})}{(1-\xi)(1-\xi^2)\dots(1-\xi^x)}$$
for $x < N$ and the $\xi$-Pochhammer symbols
$$(z;\xi)_n = (1-z)(1-z\xi)\dots(1-z\xi^{n-1}).$$
The $\xi$-Pochhammer symbol and the Gaussian binomial coefficient are related (after normalization) by the discrete Fourier transform.
To see this, recall that the $\xi$-Pochhammer symbol has the following $\xi$-binomial expansion
$$(z;\xi)_{n-1} = \sum_{x=0}^{n-1}(-z)^x\xi^{x(x-1)/2}\binom{n-1}{x}_\xi,$$
see e.g. \cite[page 11]{koekoek}.
It follows that for $n\leq N$, the discrete Fourier transform of 
$$f(x) = (-z)^x\xi^{x(x-1)/2}\binom{n-1}{x}_\xi$$
on $\mathbb{Z}/N\mathbb{Z}$ is the function
$$\wh f(x) = (z\xi^{-x};\xi)_{n-1}.$$
This leads to the following result.

\begin{comment}
\bth{starting basis}
Let $N=4m+2-d$ for some $m>0$ and $1\leq d\leq 4$.
Fix an integer $m\leq a\leq (N-1)/2$, let $r = 4(a-m)+d$, and define the function $\psi_a\in L^2(\mathbb{Z}/N\mathbb{Z})$ whose support is contained in $[-a,a]$ by $\psi_a(-a)=1$ and
\begin{align}
\psi_a(x)
&= (-1)^{(r+1)x}\xi^{x(x+1)/2}\xi^{-\frac{1}{2}(a^2+x + r(x+a))}\binom{2a-r+1}{x+a}_\xi\\\nonumber
&= (-1)^{(r+1)x}\prod_{k=1}^{x+a}\frac{\sin(\pi (2a+2-r-k)/N)}{\sin(\pi k/N)},\quad -a < x < a.
\end{align}

Then for any integer $m\leq a\leq (N-1)/2$ the space $\mathcal C_a$ has a basis of the form
$$B = \{\psi_a(x),\psi_a(x-1),\dots,\psi_a(x-r+1)\}.$$
\eth
\begin{proof}
For each integer $0\leq k < r$, the support of $\psi_a(x-k)$ is contained in $\{-a+k,a-r+k\}$, with $\psi_a(k-a)\neq0$.  It follows easily from this that the elements of $B$ are linearly independent.
Moreover, the dimension of $\mathcal C_m$ is $r$ by Proposition \ref{pCa}, so it suffices to show that $\psi_a(x),\dots,\psi_a(x-r+1)$ are all elements of $\mathcal C_m$.
Since their supports are all contained in $\{-a,1-m,m-1,a\}$, it suffices to show that the Fourier transforms have support contained within the same set.

The Fourier transform of $\psi_a(x-k)$ is 
\begin{align*}
\xi^{kx}\wh\psi_a(x)
  &  = \xi^{kx-\frac{1}{2}(a^2+ar)}((-1)^r\xi^{-x-\frac{1}{2}(1+r)};\xi)_{2a+2-r}.
\end{align*}
which has support contained in 
$[-a,a]$.
\end{proof}
\end{comment}

\bth{starting basis}
Let $N=4m+2-d$ for some $m>0$ and $1\leq d\leq 4$.
Fix an integer $m\leq a\leq (N-1)/2$, let $r = 4(a-m)+d$, and define the function $\psi_a\in L^2(\mathbb{Z}/N\mathbb{Z})$ with support $[-a,a]$ by
\begin{align*}
\psi_a(x)
&= \frac{1}{(2i)^{N-2a-1}}\xi^{-\frac{1}{4}N(N-2a-1) + ax}(\xi^{a+1-x};\xi)_{N-2a-1}\\
&=\xi^{(r-1)x/2}\prod_{k=1}^{N-2a-1}\sin\left(\frac{\pi (a+k-x)}{N}\right).
\end{align*}

Then for any integer $m\leq a\leq (N-1)/2$ the space $\mathcal C_a$ has a basis of the form
$$B := \{\psi_a(x),\xi^{-x}\psi_a(x),\dots,\xi^{(1-r)x}\psi_a(x)\}.$$
\eth
\begin{proof}
For each integer $0 \leq k < r$, the support of $\xi^{-kx}\psi_a(x)$ is $[-a,a]$.
Furthermore, its inverse Fourier transform is
$$\frac{1}{(2i)^{N-2a-1}}\xi^{-\frac{1}{4}N(N-2a-1)}(-1)^x\xi^{(a+1)x}\xi^{x(x-1)/2}\binom{N-2a-1}{x+a-k}_\xi,$$
which has support $[-a+k,a+1-r+k]$.
It follows that $\xi^{-kx}\psi_a(x)\in \mathcal C_a$ for all $0\leq k < r$.
Further, the supports of the inverse Fourier transforms imply that the $r$ elements of $B$ are all linearly independent.
The dimension of $\mathcal C_m$ is $r$ by Proposition \ref{pCa}, so $B$ is a basis.
\end{proof}

%\subsection{Functions of Extremal Support}
%\label{sec:extremal}

One version of the Uncertainty Principle for functions on $\mathbb{Z}/N\mathbb{Z}$ is to compare the relative sizes of the support of a function $f(x)$ and its Fourier transform $\wh f(x)$.
The Donoho--Stark uncertainty principle \cite{Donoho} states that 
$$\lvert \supp(f)\rvert \lvert\supp(\wh f)\rvert \geq N.$$
Gr\"unbaum obtained a lower bound of the product of the expectations of the squares of the position and momentum operators at a given state in \cite{grunbaum2}

A stronger version of the Donoho--Stark inequality is possible in the case when $N=p$ is prime.
In this setting, Tao \cite{tao} proved that
$$\lvert \supp(f)\rvert + \lvert\supp(\wh f)\rvert \geq p+1.$$
Moreover, Tao showed that this inequality is sharp and that given any two subsets $A,B\subseteq \mathbb{Z}/N\mathbb{Z}$ with $\lvert A\rvert + \lvert B \rvert = p+1$, there exists a function $f(x)$ whose support is $A$ and whose Fourier transform has support $B$.
We refer to such a function as a \vocab{function of extremal support} on $\mathbb{Z}/p\mathbb{Z}$.
We extend this definition to non-prime values of $N$: 
\bde{extremal-support}
We say that a function $f(x)$ on $\mathbb{Z}/ N \mathbb{Z}$ is a function of extremal support if for all functions $g(x)$ on $\mathbb{Z}/N\mathbb{Z}$, 
$$
\supp(g)\subseteq \supp(f) \quad \mbox{and} \quad \supp(\hat g)\subseteq\supp(\hat f)
\quad \Rightarrow \quad g(x) = \mu f(x) \; \; \mbox{for some} \; \; \mu \in \Cset.
$$
\ede
\bth{extremal-support-psi}
Let $N=4m+2-d$ for some $m>0$ and $1\leq d\leq 4$.
The elements of the basis  of $\mathcal C_m$ stated in \thref{starting basis} are all functions of extremal support.
\eth
\begin{proof}
Let $0\leq j < r$ and $f(x) =  \xi^{-jx} \psi_m(x)$.
The support of $f(x)$ is $[-m,m]$ and the support of $\hat f(x)$ is $[-m+r-1-k,m-k]$.
Suppose that $g(x)$ is a function on $\mathbb{Z}/N\mathbb{Z}$ with the property that
$$
\supp(g)\subseteq\supp(f)
\quad\text{and}\quad
\supp(\wh g)\subseteq\supp(\wh f).
$$
Then $g\in\mathcal C_m$, so we can write
$$g(x) = \sum_{k=0}^{d-1}c_k\xi^{-kx}\psi_m(x).$$
By taking the discrete Fourier transform of 
 both sides and comparing the supports, we see that $c_k=0$ for $k\neq j$, and therefore $g(x) = c_jf(x)$.
This proves that $f(x)$ is a function of extremal support.
\end{proof}

\sectionnew{Eigenvalue inequalities and strict separation of spectra}
\label{sec:separation}

As we proved above, the spectrum of the operator $J$ can be divided into two distinct pools: a collection of eigenvalues $\mu_1,\dots,\mu_r$ of $J$ of multiplicity $1$ (which are precisely the eigenvalues of $J$ restricted to $\mathcal C_a$) and a collection $\nu_1,\dots,\nu_s$ of eigenvalues of $J$ of multiplicity two (which are precisely eigenvalues of $J$ restricted to the complement $[-a,a]'$).  In this section, we prove that these pools of eigenvalues are \textit{strictly separated} in the sense that there is a real number with the property that all multiplicity two eigenvalues lie below it and all multiplicity one eigenvalues lie above.  

In order to prove this, we first consider the matrix representation of the operator $J$ restricted to $\mathcal C_a$, in terms of the extremal basis defined in the previous section.
\ble{restricted J}
Let
$$\wt A(y) = -\left[\cos\left(\frac{\pi (2y-2a+1-r)}{N}\right)+\cos\left(\frac{\pi (2a+1)}{N}\right)\right],$$
$$\wt B(y) = 2\cos\left(\frac{\pi (2y+1-r)}{N}\right),$$
and
$$\wt C(y) = -\left[\cos\left(\frac{\pi (2y+2a+1-r)}{N}\right)+\cos\left(\frac{\pi (2a+1)}{N}\right)\right].$$
Then with respect to the extremal basis $\{\xi^{-kx}\psi_a(x)\}_{k=0}^{r-1}$, the matrix $J$ has the matrix representation
\begin{equation}\label{eqn:J|C_a}
J|_{\mathcal C_a} = \left[\begin{array}{ccccc}
\wt B(0) & \wt A(0) &  0  & \dots & 0\\
\wt C(1) & \wt B(1) & \wt A(1) & \dots & 0\\
 0  & \wt C(2) & \wt B(2) & \dots & 0\\
\vdots & \vdots & \vdots & \ddots & \vdots\\
0  &  0  &  0  & \dots & \wt B(r-1)
\end{array}\right]
\end{equation}
\ele
\begin{proof}
First, note that
$$\frac{\psi_a(x+1)}{\psi_a(x)} = \xi^{(r-1)/2} \frac{\sin\left(\frac{\pi(a-x)}{N}\right)}{\sin\left(\frac{\pi(x+a+1)}{N}\right)}\quad\text{and}\quad \frac{\psi_a(x-1)}{\psi_a(x)} = \xi^{-(r-1)/2} \frac{\sin\left(\frac{\pi(x+a)}{N}\right)}{\sin\left(\frac{\pi(a+1-x)}{N}\right)}.$$
Therefore for $0\leq k < r$, 
\begin{align*}
J(\xi^{-kx}\psi_a(x))
  & = \left[\cos\left(\frac{\pi(2x+1)}{N}\right)-\cos\left(\frac{\pi(2a+1)}{N}\right)\right]\xi^{-k(x+1)}\psi_a(x+1)\\
  & - 2\cos\left(\frac{\pi(2a+1)}{N}\right)\cos\left(\frac{2\pi x}{N}\right)\xi^{-kx}\psi_a(x)\\
  & + \left[\cos\left(\frac{\pi(2x-1)}{N}\right)-\cos\left(\frac{\pi(2a+1)}{N}\right)\right]\xi^{-k(x-1)}\psi_a(x-1)\\
  & = \left[2\xi^{(r-1-2k)/2}\sin^2\left(\frac{\pi(a-x)}{N}\right) - 2\cos\left(\frac{\pi(2a+1)}{N}\right)\cos\left(\frac{2\pi x}{N}\right)\right.\\
  & \left. + 2\xi^{-(r-1-2k)/2}\sin^2\left(\frac{\pi(a+x)}{N}\right)\right]\xi^{-kx}\psi_a(x)\\
  & = -\left[\cos\left(\frac{\pi(r-1-2k+2a)}{N}\right)+\cos\left(\frac{\pi(2a+1)}{N}\right)\right]\xi^{-(k+1)x}\psi_a(x)\\
  & + 2\cos\left(\frac{\pi(r-1-2k)}{N}\right)\xi^{-kx}\psi_a(x)\\
  & - \left[\cos\left(\frac{\pi(r-1-2k-2a)}{N}\right)+\cos\left(\frac{\pi(2a+1)}{N}\right)\right]\xi^{-(k-1)x}\psi_a(x)\\
  & = -\left[\cos\left(\frac{\pi(r-1-2k+2a)}{N}\right)+\cos\left(\frac{\pi(2a+1)}{N}\right)\right]\xi^{-(k+1)x}\psi_a(x)\\
  & + 2\cos\left(\frac{\pi(r-1-2k)}{N}\right)\xi^{-kx}\psi_a(x)\\
  & - \left[\cos\left(\frac{\pi(r-1-2k-2a)}{N}\right)+\cos\left(\frac{\pi(2a+1)}{N}\right)\right]\xi^{-(k-1)x}\psi_a(x)
\end{align*}
\end{proof}
Using the previous lemma, we can obtain upper and lower bounds for the simple eigenvalues of $J$.
\bpr{simple eigenbound}
Let $a$ be an integer with $(N-1)/4 \leq a < (N-1)/2$.  Then the simple eigenvalues $\mu_1,\dots,\mu_r$ of $J$ satisfy
$$2\cos\left(\frac{\pi}{N}\right)+4\cos\left(\frac{\pi(2a+1)}{N}\right)\leq \mu_j\leq 2\cos\left(\frac{\pi}{N}\right)-4\cos\left(\frac{\pi(2a+1)}{N}\right)$$
for all $1\leq j\leq r$.
\epr
\begin{proof}
By conjugating the matrix representation of $J|_{\mathcal C_a}$ by an appropriate diagonal matrix, we can symmetrize the tridiagonal matrix, obtaining the matrix
$$\left[\begin{array}{ccccc}
\wt B(0) & \sqrt{\wt A(0)\wt C(1)} &  0  & \dots & 0\\
\sqrt{\wt A(0)\wt C(1)} & \wt B(1) & \sqrt{\wt A(1)\wt C(2)} & \dots & 0\\
 0  & \sqrt{\wt A(1)\wt C(2)} & \wt B(2) & \dots & 0\\
\vdots & \vdots & \vdots & \ddots & \vdots\\
0  &  0  &  0  & \dots & \wt B(r-1)
\end{array}\right].$$
In particular, the eigenvalues of this matrix are $\mu_1,\dots,\mu_r$.
The Gershgorin's Circle Theorem \cite{Viazovska} implies that for each $1\leq j\leq r$ there exists a $0\leq k < r$ with
$$|\mu_j-\wt B(k)| \leq \sqrt{\wt A(k-1)\wt C(k)} + \sqrt{\wt A(k)\wt C(k+1)}.$$
Note that we can factor
$$\wt A(k)\wt C(k+1) = 
\left(1+\cos\left(\frac{\pi(2y+2-r)}{N}\right)\right)\left(\cos\left(\frac{\pi(2y+2-r)}{N}\right)-\cos\left(\frac{\pi r}{N}\right)\right).$$
Therefore 
\begin{align*}
\left[\cos\left(\frac{\pi(2y+2-r)}{N}\right)-1-2\cos\left(\frac{\pi(2a+1)}{N}\right)\right]^2-\wt A(k-1)\wt C(k)\\
= 2\left(1+\cos\left(\frac{\pi(2a+1)}{N}\right)\right)^2\left(1-\cos\left(\frac{\pi(2y+2-r)}{N}\right)\right)\geq 0.
\end{align*}
Since the quantity inside the square is positive, this shows
$$\sqrt{\wt A(k)\wt C(k+1)}\leq \cos\left(\frac{\pi(2y+2-r)}{N}\right)-1-2\cos\left(\frac{\pi(2a+1)}{N}\right).$$
Consequently,
$$\sqrt{\wt A(k-1)\wt C(k)}\leq \cos\left(\frac{\pi(2y-r)}{N}\right)-1-2\cos\left(\frac{\pi(2a+1)}{N}\right).$$
Putting this together with the inequality from Gershgorin's Circle Theorem, we obtain the spectral lower bound
\begin{align*}
\mu_j
 & \geq \min_{0\leq k < r} \left(\wt B(k) - \sqrt{\wt A(k-1)\wt C(k)} - \sqrt{\wt A(k)\wt C(k+1)}\right)\\
 & \geq \min_{0\leq k < r} \cos\left(\frac{\pi(2y+1-r)}{N}\right)\left(2-2\cos\left(\frac{\pi}{N}\right)\right)+2+4\cos\left(\frac{\pi(2a+1)}{N}\right)\\
 & \geq 2\cos\left(\frac{\pi}{N}\right) + 4\cos\left(\frac{\pi(2a+1)}{N}\right).
\end{align*}
Likewise, we obtain the spectral upper bound
\begin{align*}
\mu_j
 & \leq \max_{0\leq k < r} \left(\wt B(k) + \sqrt{\wt A(k-1)\wt C(k)} + \sqrt{\wt A(k)\wt C(k+1)}\right)\\
 & \leq \max_{0\leq k < r} \cos\left(\frac{\pi(2y+1-r)}{N}\right)\left(2+2\cos\left(\frac{\pi}{N}\right)\right)-2-4\cos\left(\frac{\pi(2a+1)}{N}\right)\\
 & \leq  2\cos\left(\frac{\pi}{N}\right) - 4\cos\left(\frac{\pi(2a+1)}{N}\right).
\end{align*}
\end{proof}

Gershgorin's Circle Theorem can also be used to obtain bounds for the double eigenvalues of $J$.
\bpr{double eigenbound}
Let $a$ be an integer with $(N-1)/4 \leq a < (N-1)/2$.  Then the double eigenvalues $\nu_1,\dots,\nu_r$ of $J$ satisfy
$$-2\cos\left(\frac{\pi}{N}\right)+4\cos\left(\frac{\pi(2a+1)}{N}\right)\leq \nu_j\leq 2\cos\left(\frac{\pi}{N}\right)+4\cos\left(\frac{\pi(2a+1)}{N}\right)$$
for all $1\leq j\leq r$.
\epr
\begin{proof}
Let $A(x)$ and $B(x)$ be the functions in the definition of the operator $J$, ie.
\begin{align*}
  A(x) &= \cos\left(\frac{\pi (2x+1)}{N}\right)-\cos\left(\frac{\pi (2a+1)}{N}\right),\\
  B(x) &= - 2\cos\left(\frac{\pi (2a+1)}{N}\right)\cos\left(\frac{2\pi x}{N}\right).
\end{align*}
Then the transformation $J$ restricted to the complement $[-a,a]'$ of $[-a,a]$ is represented in terms of the standard basis by the matrix
$$\left[\begin{array}{ccccc}
B(a+1) & A(a+1) &  0  & \dots & 0\\
A(a+1) & B(a+2) & A(a+2) & \dots & 0\\
 0  & A(a+2) & B(a+3) & \dots & 0\\
\vdots & \vdots & \vdots & \ddots & \vdots\\
0  &  0  &  0  & \dots & B(N-a-1)
\end{array}\right].
$$
The eigenvalues of this matrix are exactly the eigenvalues $\nu_1,\dots,\nu_s$.
Applying Gershgorin's Circle Theorem, we get that for all $1\leq j\leq s$ there must exist $a<k<N-a$ with
$$|\nu_j-B(k)|\leq |A(k-1)|+|A(k)|.$$
Furthermore for $a<k<N-a$
$$|A(k-1)|+|A(k)| = -2\cos\left(\frac{2\pi k}{N}\right)\cos\left(\frac{\pi}{N}\right)+2\cos\left(\frac{\pi(2a+1)}{N}\right).$$
Therefore
\begin{align*}
\nu_j
 & \geq \min_{a < k < N-a} \left(B(k)-|A(k-1)|-|A(k)|\right)\\
 & = \min_{a < k < N-a} \cos\left(\frac{2\pi k}{N}\right)\left[2\cos\left(\frac{\pi}{N}\right)-2\cos\left(\frac{\pi(2a+1)}{N}\right)\right] - 2\cos\left(\frac{\pi(2a+1)}{N}\right)\\
 & \geq -2\cos\left(\frac{\pi}{N}\right) + 4\cos\left(\frac{\pi(2a+1)}{N}\right).
\end{align*}
and also
\begin{align*}
\nu_j
 & \leq \max_{a < k < N-a} \left(B(k)+|A(k-1)|+|A(k)|\right)\\
 & = \max_{a < k < N-a} \cos\left(\frac{2\pi k}{N}\right)\left[-2\cos\left(\frac{\pi}{N}\right)-2\cos\left(\frac{\pi(2a+1)}{N}\right)\right] + 2\cos\left(\frac{\pi(2a+1)}{N}\right)\\
 & \leq 2\cos\left(\frac{\pi}{N}\right) + 4\cos\left(\frac{\pi(2a+1)}{N}\right).
\end{align*}
\end{proof}

\begin{figure}
    \centering
    \includegraphics[width=0.8\linewidth]{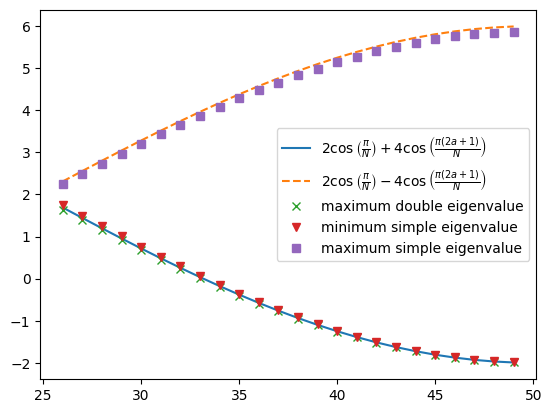}
    \caption{Plot of the largest double eigenvalue and the largest and smallest simple eigenvalues of $J$ for various values of $a$ when with $N=101$, along with the boundaries for the maximum and minimum of the simple spectra.  As one observed, the actual largest and smallest values closely resemble the computed bounds.}
    \label{fig:eigengraph}
\end{figure}

When combined, the previous propositions show that whether an eigenvalue is simple or double is determined by whether it is above or below $2+4\cos\left(\frac{\pi (2a+1)}{N}\right)$.  In particular, this proves part (e) of the Spectral Theorem.  In practice, the lower bound for the double eigenvalues we obtained is not very good.  However, the upper bound of the double eigenvalues, and the lower and upper bounds of the simple eigenvalues turn out to be remarkably close to the numerical values of the minimum and maximum eigenvalues, especially for large values of $N$ and $a$.
The graph in Figure \ref{fig:eigengraph} shows the extremal eigenvalues of $J$ for $N=101$ as $a$ ranges between $(N-1)/4$ and $(N-1)/2$, along with the theoretical bounds established in this section.  As one sees, the curves themselves are very similar in nature.

\sectionnew{The cases of $\mathcal C_a$ of dimensions $1$, $2$, $3$ and $4$}
\label{sec:low dim}
The extremal basis for $\mathcal C_a$ found in Theorem \ref{tstarting basis} allows us to obtain explicit expressions for the $F_N$ and $J$ joint eigenbases of those spaces in the cases when 
$$
1 \leq \dim \mathcal C_a \leq 4.
$$
The essential reason is that if one diagonalizes the matrix representing the restriction $J|_{\mathcal C_a}$ in Equation \eqref{eqn:J|C_a}, then one obtains an explicit expression for the eigenvectors in the joint eigenbasis of $F_N$ and $J$ in terms of the extremal basis for $\mathcal C_a$.  When $\dim\mathcal C_a$ is small enough, this process is clean enough to be carried out by hand.
These results also provide the bases cases for the inductive proof of \prref{eigenvalues of restriction}. 

Denote once again $N=4m+2-d$ with $1\leq d\leq 4$ and consider the case $a=m$, so $\dim \mathcal C_m = d$.
\begin{enumerate}[\textbf{Case} 1:]
\item $(N=4m+1)$.  In this case $\mathcal C_m$ is one-dimensional, so 
\begin{equation}\label{case 1}
\psi_m(x) = \prod_{k=1}^{2m}\sin\left(\frac{\pi (m+k-x)}{N}\right)
\end{equation}
is already an eigenfunction.  The corresponding eigenvalue is $\sqrt{N}$.
\item $(N=4m)$.  In this case $\mathcal C_m$ is two dimensional and spanned by $\psi_m(x)$ and $\xi^{-x}\psi_m(x)$.  In particular, up to constant multiples it contains unique even and odd functions, given by 
\begin{equation}\label{case 2a}
\frac{1}{2}(1 + \xi^{-x})\psi_m(x) = \cos\left(\frac{\pi x}{N}\right)\prod_{k=1}^{2m-1}\sin\left(\frac{\pi (m+k-x)}{N}\right)
\end{equation}
and
\begin{equation}\label{case 2b}
\frac{1}{2i}(1 - \xi^{-x})\psi_m(x) = \sin\left(\frac{\pi x}{N}\right)\prod_{k=1}^{2m-1}\sin\left(\frac{\pi (m+k-x)}{N}\right),
\end{equation}
respectively.  Eigenfunctions of $F_N$ with real eigenvalues are even, while those with imaginary eigenvalues are odd, so these must each be eigenfunctions of $F_N$.  The corresponding eigenvalues are $\sqrt{N}$ and $-\sqrt{N}i$.
\item $(N=4m-1)$.  In this case $\mathcal C_m$ is three dimensional and has a unique odd function (up to a constant multiple), given by 
\begin{equation}\label{case 3}
\frac{1}{2i}(1-\xi^{-2x})\psi_m(x) = \sin\left(\frac{2\pi x}{N}\right)\prod_{k=1}^{2m-2}\sin\left(\frac{\pi (m+k-x)}{N}\right),
\end{equation}
which must be an eigenfunction of $F_N$.  The corresponding eigenvalue is $-i\sqrt{N}$.
In particular, this implies
$$\wh\psi_m (x) - \wh\psi_m(x+2) = -i\sqrt{N}(1-\xi^{-2x})\psi_m(x).$$
The remaining eigenfunctions will be even with eigenvalues $\pm\sqrt{N}$, and therefore they will have to be scalar multiples of functions of the form
$$(1 + c\xi^{-x} + \xi^{-2x})\psi_m(x),$$
for some constant $c$.
Taking the discrete Fourier transform, this gives
$$\wh \psi_m(x) + c\wh\psi_m(x+1) + \wh\psi_m(x+2) = \pm\sqrt{N}(1 + c\xi^{-x} + \xi^{-2x})\psi_m(x).$$
Consequently,
{\small
$$(\wh \psi_m(x) + c\wh\psi_m(x+1) + \wh\psi_m(x+2)) = \mp i\frac{1 + c\xi^{-x} + \xi^{-2x}}{(1-\xi^{-2x})}(\wh\psi_m (x) - \wh\psi_m(x+2)).$$
}
Noting that $\wh\psi(m)\neq 0$, $\wh \psi(m+1) = 0$ and $\wh\psi(m+2) = 0$, if we evaluate this expression at $x=m$, we find
$$c = -2\cos(2\pi m/N)\pm 2\sin(2\pi m/N).$$
Thus, the eigenfunctions corresponding to $\pm\sqrt{N}$ can be taken to be
\begin{align}\label{case 3b}
&(1 + c\xi^{-x} + \xi^{-2x})\psi_m(x)\\\nonumber
&=\left(\cos\left(\frac{2\pi x}{N}\right) - \cos\left(\frac{2\pi m}{N}\right) \pm \sin\left(\frac{2\pi m}{N}\right)\right)\prod_{k=1}^{2m-2}\sin\left(\frac{\pi (m+k-x)}{N}\right).
\end{align}

\item $(N=4m-2)$.  In this case $\mathcal C_m$ is four dimensional and the eigenfunctions appear as either even functions
$$(1 + c\xi^{-x} + c\xi^{-2x} + \xi^{-3x})\psi_m(x),$$
or odd functions
$$(1 + c\xi^{-x} - c\xi^{-2x} - \xi^{-3x})\psi_m(x)$$
for some specific values of $c$.
If we have the same eigenvalue repeated in the eigenspace, then by taking their difference, we obtain an element of the space $\mathcal C_{m-1}$.
Since this space is trivial, each possible eigenvalue of $F$ occurs exactly one time.
Therefore we have two even eigenfunctions with eigenvalues $\pm\sqrt{N}$, and two odd eigenfunctions with eigenvalues $\pm i\sqrt{N}$.

To figure out the exact value of the unknown constant, we can use the fact that we have an explicit expression for the Fourier transform of $\psi_m(x)$, namely
$$\wh\psi_m(x) = (-2i)^{2m-3}\xi^{-\frac{1}{4}(4m-2)(2m-3)}(-1)^x\xi^{-(m+1)x}\xi^{-x(-x+1)/2}\binom{2m-4}{-x+m}_\xi.$$
Here we are using that $F_N^2 f(x) = f(-x)$ for any function $f(x)$.
In particular, $\wh\psi_m(m) \neq 0$ and $\wh\psi_m(m+k) = 0$ for $0<k<4$, so the eigenfunctions may be obtained by evaluation at the point $x=m$.
Specifically, for the eigenvalues $\pm\sqrt{N}$, we must have
$$\wh\psi_m(m) = \pm\sqrt{N}(1 + c\xi^{-m} + c\xi^{-2m} + \xi^{-3m})\psi_m(m),$$
which says
$$c = \frac{\wh\psi_m(m)}{\pm2\sqrt{N}\cos(\pi m/N)\xi^{-3m/2}\psi_m(m)}-\frac{\cos(3\pi m/N)}{\cos(\pi m/N)}.$$
A similar expression holds for the odd eigenfunction expression with the eigenvalues $\pm i\sqrt{N}$.
In particular, we have the four eigenfunctions
\begin{align}\label{case 4}
\left(\cos\left(\frac{3\pi x}{N}\right) + \left(\frac{\wh\psi_m(m)}{\pm\sqrt{N}\xi^{-3m/2}\psi_m(m)}-\cos\left(\frac{3\pi m}{N}\right)\right)\frac{\cos\left(\frac{\pi x}{N}\right)}{\cos\left(\frac{\pi m}{N}\right)}\right)\xi^{-3x/2}\psi_m(x),\\
\left(\sin\left(\frac{3\pi x}{N}\right) + \left(\frac{\wh\psi_m(m)}{\mp\sqrt{N}\xi^{-3m/2}\psi_m(m)}-\sin\left(\frac{3\pi m}{N}\right)\right)\frac{\sin\left(\frac{\pi x}{N}\right)}{\sin\left(\frac{\pi m}{N}\right)}\right)\xi^{-3x/2}\psi_m(x)
\end{align}
with eigenvalues $\pm \sqrt{N}$ and $\pm i\sqrt{N}$, respectively.
\end{enumerate}
\bre{prev-results}
Each of the above cases give us very concrete expressions for some of the eigenfunctions of the discrete Fourier transform $F_N$.
The single eigenfunction $\psi_m(x)$ belonging to $\mathcal C_m$ in the case $N=4m+1$ shown in Eq. \eqref{case 1} was also found (up to a constant multiple) by Kong in \cite{kong}, but expressed in the algebraically equivalent form
$$\prod_{k=m+1}^{2m} \left[\cos\left(\frac{2\pi}{N}x\right) - \cos\left(\frac{2\pi}{N}k\right)\right].$$
Kong also obtained the following similar expression for the odd eigenfunction in $\mathcal C_m$ in the case $N=4m$ appearing in Equation \eqref{case 2a}:
$$\sin\left(\frac{2\pi}{N}x\right) \prod_{k=m+1}^{2m-1}\left[\cos\left(\frac{2\pi}{N}x\right)- \cos\left(\frac{2\pi}{N}k\right)\right].
$$
However, the remaining explicit expressions in Eqs. \eqref{case 2b}, \eqref{case 3}, \eqref{case 3b}, and \eqref{case 4} are new.
Thus our expressions for the extremal eigenfunctions provide a nice extension of this collection of known results.
\ere
\sectionnew{Reconstruction for Fourier uniqueness pairs}
A \vocab{Fourier uniqueness pair} for $\mathbb{Z}/N\mathbb{Z}$ is a pair of subsets $A,B\subseteq \mathbb{Z}/N\mathbb{Z}$ with the property that $f\in L^2(\mathbb{Z}/N\mathbb{Z})$ is uniquely determined by knowing its value on $A$, along with the value of $\wh f$ on $B$. 
The study of discrete subsets of the real line which form Fourier uniqueness pairs has been a topic of many recent papers, including \cite{Ramos,Radchenko}.
One important related question is, given a Fourier uniqueness pair, how to obtain an explicit interpolation formula allowing for the reconstruction of the function.

It follows from the Spectral Theorem that, for $a\geq (N-2)/4$, the pair of identical discrete intervals $A = [-a,a]$ and $B = [-a,a]$ forms an analog to a Fourier uniqueness pair in the setting of the finite Fourier transform.
The goal of this section is to obtain a corresponding interpolation formula, 
which proves parts (a) and (b) of the Interpolation Theorem.
Put another way, we want to reconstruct the value of $f(x)$ on $[-a,a]'$ from knowing the value of $f(x)$ and $\wh f(x)$ on $[-a,a]$.
The eigendata of $J$ provides a novel solution to this problem.

As we showed in the Spectral Theorem, the eigenvalues of $J$ are
$$\mu_1,\dots,\mu_r,\nu_1,\dots,\nu_s
$$
with $r=4a+2-N$ and $s=N-2a-1$, where each $\mu_k$ is an eigenvalue of $J$ with multiplicity $1$ and each $\nu_j$ has multiplicity $2$.
Furthermore, we can choose a basis for the $\nu_j$-eigenspace of $J$ consisting of two functions $\varphi_j(x)$ and $\wt\varphi_j(x)$, supported on $[-a,a]$ and $[-a,a]'$, respectively.

The fact that $F_N$ commutes with $J$ means that $F_N$ will preserve the eigenspace of $\nu_j$.
The next proposition describes the action of $F_N$ on this space (we use the notation from the Spectral Theorem in the introduction.)
\bpr{FN-square}
Let $(N-2)/4\leq a \leq (N-1)/2$, $r:=4a+2-N$ and $s:=N-2a-1$.
For each $1\leq j < s$, there exist nonzero numbers $\alpha_j$ and $\beta_j$, either both real or both imaginary, with 
$|\alpha_j|^2+|\beta_j|^2=1$ and
$$\left[\begin{array}{c}
F_N \varphi_j(x)\\ F_N\wt\varphi_j(x)
\end{array}\right] = 
\left[\begin{array}{cc}
\alpha_j & \beta_j\\
\beta_j & -\alpha_j
\end{array}\right]
\left[\begin{array}{c}
 \varphi_j(x)\\ \wt\varphi_j(x)
\end{array}\right]
.$$
\epr
\begin{proof}
The Fourier transform acts as a unitary operator on the eigenspace of $J$ with eigenvalue $\nu_j$.
Therefore there exist complex numbers $\alpha_j$, $\beta_j$, and $\gamma_j$, with $|\alpha_j|^2+|\beta_j|^2=1$, $|\gamma_j| = 1$, and 
$$\left[\begin{array}{c}
F_N \varphi_j(x)\\ F_N\wt\varphi_j(x)
\end{array}\right] = 
\left[\begin{array}{cc}
\alpha_j & \beta_j\\
-\gamma_j\overline\beta_j & \gamma_j\overline\alpha_j
\end{array}\right]
\left[\begin{array}{c}
 \varphi_j(x)\\ \wt\varphi_j(x)
\end{array}\right].$$
Moreover, $\varphi_j(x)\notin\mathcal C_a$, so $\beta_j\neq 0$.
Since $\varphi_j(x)$ is real we know that
$$\wh{\varphi_j}(-x) = \overline{\wh\varphi_j}(x),$$
so that
$$\alpha_j\varphi_j(-x) + \beta_j\wt \varphi_j(-x) = \overline\alpha_j\varphi_j(x) + \overline\beta_j\wt \varphi_j(x).$$
By the symmetry of $J$ and the multiplicity of the eigenvalue $\nu_j$, the eigenfunctions $\varphi_j(x)$ and $\wt\varphi_j(x)$ will be either even or odd.
Furthermore, the discrete Fourier transform of a real function sends even functions to real and imaginary ones, so both $\varphi_j(x)$ and $\wt\varphi_j(x)$ have the same parity.
Since $\wt\varphi_j(x)$ is also real, this means $\alpha_j$ and $\beta_j$ are both simultaneously either purely real or purely imaginary.

Since $F_N^2=\pm I$, we see $\gamma_j=\mp 1$ and $\gamma_j\overline\alpha_j = -\alpha_j$.
If $\alpha_j$ is real, this implies that $\gamma_j=-1$ and $-\gamma_j\overline\beta_j = \beta_j$.
Likewise, if $\alpha_j$ is imaginary, then $\gamma_j=1$ and $-\gamma_j\overline\beta_j =\beta_j$.
This completes the proof.
\end{proof}
\begin{proof}[Proof of part (b) of the Interpolation Theorem] The time-band limited discrete Fourier transform $F_N^a$ commutes with the operator $J$ because $J$ commutes with $F_N$ and the projection $P_a$ (\prref{JFN-comm} and proof of \prref{repeated eigenvalues of J}). The functions $\rho_k(x)$, $1 \leq k \leq r$ are eigenfunctions of $F_N^a$ with the same eigenvalues as $F_N$ because of \eqref{tbF-F}. Since
$\wt\varphi_j(x) \in L^2([-a,a]'),$
$$
F_N^a (\wt\varphi_j(x)) = 0, \quad \forall 1 \leq j \leq s.
$$

Part (b) of the Interpolation Theorem
and the property that $\varphi_j(x) \in L^2([-a,a])$
imply
$$
F_N^a \varphi_j(x) = 
P_a F_N P_a \varphi_j(x) = P_a F_N \varphi_j(x) 
= P_a ( \alpha_j \varphi_j(x) + \beta_j \wt\varphi_j(x)) = \alpha_j \varphi_j(x)
$$
for all $1 \leq j \leq s$.
Lastly, since $P_a\rho_k=\rho_k$, $P_a\varphi_j=\varphi_j$ and $P_a\wt\varphi_j=0$, we calculate
$$B_N^a\rho_k = P_aF_N^*P_aF_NP_a\rho_k = P_aF_N^*F_N\rho_k = \overline\nu_k\nu_k P_a\rho_k = N^2\rho_k,$$
$$B_N^a\wt\varphi_j=P_aF_N^*P_aF_NP_a\wt\varphi_j = 0,$$
and also
$$B_N^a\varphi_j=P_aF_N^*P_aF_NP_a\varphi_j = P_aF_N^*P_aF_N\varphi_k=\alpha_jP_aF_N^*\varphi_j=\lvert\alpha_j\rvert^2\varphi_j.$$
This completes the proof of this part of the Interpolation Theorem.
\end{proof}

\begin{proof}[Proof of part (c) of the Interpolation Theorem]
For any $f(x)\in L^2(\mathbb{Z}/N\mathbb{Z})$, we can expand
$$f(x) = \sum_{j=1}^{r} b_k\rho_k(x) + \sum_{k=j}^s (c_j\varphi_j(x) + d_j\wt\varphi_j(x)).$$
Alternatively, we can expand the Fourier transform
$$\wh f(x) = \sum_{k=1}^{r} \wt b_k\rho_k(x) + \sum_{j=1}^s(\wt c_j\varphi_j(x) + \wt d_j\wt\varphi_j(x)).$$

Comparing coefficients, we see that
$$\left[\begin{array}{c}
\wt c_j\\\wt d_j
\end{array}\right] = 
\left[\begin{array}{cc}
\alpha_j & \beta_j\\
\beta_j  & -\alpha_j
\end{array}\right]
\left[\begin{array}{c}
 c_j\\ d_j
\end{array}\right]
.$$
In particular,
$$d_j = \frac{1}{\beta_j}(\wt c_j - \alpha_j c_j),$$
or equivalently
$$d_j = \frac{1}{\beta_j}\sum_{y\in [-a,a]}(\wh f(y) - \alpha_j f(y))\varphi_j(y).$$
Therefore
$$f(x)  = \sum_{k=1}^s\wt \varphi_j(x)\left(\frac{1}{\beta_j}\sum_{y\in [-a,a]} (\wh f(y) - \alpha_j f(y))\varphi_j(y)\right)\quad\text{for all}\ x\notin[-a,a].$$
If we replace $f(x)$ with $\wh f(x)$, we also get 
$$\wh f(x)  = \sum_{j=1}^s\wt \varphi_j(x)\left(\frac{1}{\beta_j}\sum_{y\in [-a,a]} (Nf(-y) - \alpha_j \wh f(y))\varphi_j(y)\right)\quad\text{for all}\ x\notin[-a,a].$$

This defines the values of $f(x)$ outside the discrete interval $[-a,a]$, using only the values of $f(x)$ and $\wh f(x)$ inside the interval.
\end{proof}
\sectionnew{Relation to theta functions}
\label{sec:theta-functions}
The eigenfunctions of $J$ have an interesting geometric connection.
To see this, consider the theta function $\theta: \mathbb{C}\times \mathbb{H} \rightarrow\mathbb{C}$ defined by
$$\theta(x,\tau) = \sum_{n\in\mathbb{Z}}\exp(i\pi \tau(x+nN)^2/N),$$
where $\mathbb{H}$ denotes the upper half plane.
In the special case $N=2$, the functions $\theta(0,\tau) = \vartheta_{00}(0,2\tau)$ and $\theta(1,\tau) = \vartheta_{10}(0,2\tau)$ are theta constants.
In general,
$$\theta(x,\tau) = e^{i\pi \tau x^2/N}\vartheta(x\tau;N\tau)$$
for
$$\vartheta(z,\tau) = \sum_{n\in\mathbb{Z}}\exp(i\tau n^2+2\pi i nz)$$
the Jacobi theta function.
Jacobi's identity says
$$\vartheta(z/\tau,-1/\tau) = e^{i\pi z^2/\tau}\sqrt{-i\tau}\vartheta(z,\tau).$$
Therefore,
%$$\theta(x,-1/\tau) = e^{-i\pi x^2/\tau N}\vartheta((-x/N)/(\tau/N),-1/(\tau/N)) = \sqrt{-i\tau/N}\vartheta(-x/N,\tau/N)$$
$$\theta(x,-1/\tau) = \sqrt{-i\tau/N}\vartheta(-x/N,\tau/N).$$

As expected, $\theta(x,\tau)$ has several nice other algebraic properties.
For example $\theta(\tau,x)$ is periodic in both of its variables with period $N$, i.e.,
$$\theta(x+N,\tau) = \theta(x,\tau)\quad\text{and}\quad\theta(x,\tau+N)=\theta(x,\tau).$$

There is also a very nice relationship between $\theta(\tau,x)$ and its finite Fourier transform in the variable $x$.
\ble{DFT-theta}
Viewed as a function on $\mathbb{Z}/N\mathbb{Z}$, the discrete Fourier transform of $\theta(x,\tau)$ is given by
$$\wh\theta(x,\tau) = \sqrt{\frac{N}{-i\tau }}\theta(x,-1/\tau) = \vartheta(-x/N,\tau/N)\quad\text{for all} \; \;  x\in\mathbb{Z}/N\mathbb{Z}.$$
\ele
\begin{proof}
Recall the formula for the Fourier transform of a complex Gaussian:
$$\int_{\mathbb{R}} e^{i\pi \tau x^2/N}e^{-2\pi ikx/N}dx = \sqrt{\frac{N}{-i\tau}} e^{-i\pi k^2/\tau N}.$$
Therefore, for any integer $x$, we divide up the integral to find
$$\int_0^N\theta(x,\tau) e^{-2\pi i kx/N}dx = \sqrt{\frac{N}{-i\tau }} e^{-i\pi k^2/\tau N},$$
so that
$$\theta(x,\tau) = \frac{1}{\sqrt{-iN\tau }} \sum_{k=0}^{N-1}\theta(x,-1/\tau)e^{2\pi i kx/N}.$$
This completes the proof.    
\end{proof}
\begin{proof}[Proof of part (d) of the Interpolation Theorem]
The Interpolation Formula tells us that for any function $f(x)\in L^2(\mathbb{Z}/N\mathbb{Z})$ we can write
$$f(x) = \sum_{y\in [-a,a]}v_y(x)f(y) + w_y(x)\wh f(y),\quad\text{for all}\ x\notin [-a,a],$$
where 
$$
v_y(x) = \sum_{k=1}^s\frac{-\alpha_k}{\beta_k}\varphi_k(y)\wt\varphi_k(x)
\quad\text{and}\quad
w_y(x) = \sum_{k=1}^s\frac{1}{\beta_k}\varphi_k(y)\wt\varphi_k(x).
$$

Then for all $\tau\in\mathbb{H}$ and all $x\in[-a,a]'$,
$$\theta(x,\tau) = \sum_{y\in [-a,a]}\left( v_y(x)\theta(y,\tau)+\sqrt{\frac{N}{-i\tau}}w_y(x)\theta(y,-1/\tau)\right).$$
This allows us to express $v_y(x)$ and $w_y(x)$ in terms of Wronskians of theta functions, on the interval $[-a,a]'$.
Specifically, we can write for $y\in[-a,a]$ and $x\in [-a,a]'$
$$
v_y(x) = \frac{W\left(\theta(-a,\tau),\dots \theta(x,\tau),\dots,\theta(a,\tau),\vartheta(-a/N,\tau/N),\dots,\vartheta(a/N,\tau/N)\right)}{W\left(\theta(-a,\tau),\dots,\theta(a,\tau),\vartheta(-a/N,\tau/N),\dots,\vartheta(a/N,\tau/N)\right)},
$$
where $\theta(x,\tau)$ is occuring in the $(y+a+1)$'th position, and 
$$
w_y(x) = \frac{W\left(\theta(-a,\tau),\dots,\theta(a,\tau),\vartheta(-a/N,\tau/N),\dots,\theta(x,\tau),\dots,\vartheta(a/N,\tau/N)\right)}{W\left(\theta(-a,\tau),\dots,\theta(a,\tau),\vartheta(-a/N,\tau/N),\dots,\vartheta(a/N,\tau/N)\right)},  
$$
where $\theta(x,\tau)$ is occuring in the $(y+3a+2)$'th position.
In particular, these two Wronskian expressions are constant in the value $\tau$.
\end{proof} 
\bex{N2} Consider the interesting special case of 
$$
N=2 \quad \mbox{and} \quad a=0,
$$
where one can work out the values of $v_0(1)$ and $w_0(1)$ directly from the basic definitions, instead of from the formula above.
In this case the interpolation formula can be used to obtain properties of theta functions.

For every function $f: \mathbb{Z}/2 \mathbb{Z} \rightarrow \mathbb{C}$, we have
$$f(1) = v_0(1)f(0) + w_0(1)\wh f(0).$$
Since $\wh f(0) = f(0)+f(1)$, this says
$$f(1) = v_0(1)f(0) + w_0(1)f(0) + w_0(1)f(1).$$
Therefore, $v_0(1)=-1$ and $w_0(1)=1$.

Adopting a standard notation, we will write $\theta_2(\tau)=\vartheta_{10}(0,\tau)$ and $\theta_3(\tau)=\vartheta_{00}(0,\tau)$.
Then for $N=2$, we have $\theta(0,\tau)=\theta_3(2\tau)$ and $\theta(1,\tau)=\theta_2(2\tau)$.
Leveraging Jacobi's identity, the Wronskian expressions above therefore say
$$\frac{W(\theta_2(2\tau),\theta_3(\tau/2))}{W(\theta_3(2\tau),\theta_3(\tau/2))} = -1
\quad\text{and}\quad
\frac{W(\theta_3(2\tau),\theta_2(\tau/2))}{W(\theta_3(2\tau),\theta_3(\tau/2))}=1.$$
When the left hand sides are expnded, one obtains the identities 
$$
2(\theta_2'(2\tau)+\theta_3'(2\tau))'\theta_3(\tau/2)-\frac{1}{2}(\theta_2(2\tau)+\theta_3(2\tau))\theta_3'(\tau/2)
=0
$$
and
$$
2\theta_3'(2\tau)(\theta_2(\tau/2)-\theta_3(\tau/2))-\frac{1}{2}\theta_3(2\tau)(\theta_2'(\tau/2)-\theta_3'(\tau/2)) = 0.
$$
They can be derived from the formulas
$$\theta_3(\tau/2) = \theta_3(2\tau) + \theta_2(2\tau)
\quad\text{and}\quad
\theta_3(2\tau) = \theta_2(\tau/2)-\theta_3(\tau/2),$$
which in turn are consequences of Landen's transformation equations (see \cite[pp. 20, Exercise 2]{Lawden}). 
\eex

%%%%%%%%%%%%%%%%%%%%%%%%%%%%%%%%%%%%%%%%%%%%%%%%%%%%%%%%%%%%%%%%%%%%%%%%%%%%%%
\end{document}